\documentclass[12pt,twoside]{amsart}
\usepackage{amsmath}
\usepackage{amssymb}
\usepackage{amscd}
\usepackage{xypic}
\xyoption{all}
\title[Serre Intersection Multiplicity Conjecture and Hodge theory]{Serre Intersection Multiplicity Conjecture and Hodge theory}

\author{Mohammad Reza Rahmati}
\thanks{}
\address{ABDUS SALAM SCHOOL OF MATHEMATICAL SCIENCES, Pakistan
\hfill\break 
\hfill\break \\
\hfill\break }
\email{mrahmati@cimat.mx, rahmati@sms.edu.pk}

\newcommand{\comments}[1]{}

\newtheorem{theorem}{Theorem}[section]
\newtheorem{proposition}[theorem]{Proposition}



\begin{document}

\begin{abstract}
We explain intersection multiplicity defined by J. P. Serre, in terms of the Poincare product in Hodge theory by a modification of the chern character map. We also discuss a formulation of the Euler characteristic via the action of correspondences on the Chow groups of projective varieties, assuming the Grothendieck Standard conjectures over $\mathbb{Q}_l$.
\end{abstract}

\maketitle


\section*{Introduction}
\vspace{0.3cm}

Let $M,N$ be finitely generated modules over a regular local ring $A$ such that $M \otimes_A N$ has finite length. J. P. Serre \cite{S}, defines their intersection multiplicity as 

\begin{equation}
\chi^A(M,N) := \displaystyle{\sum (-1)^i l(\text{Tor}_i^A(M,N)})
\end{equation}

\noindent
He proves the basic fact that in this case $\dim M + \dim N \leq \dim A$, will hold and makes the following question, known as Serre Multiplicity conjecture.

\vspace{0.3cm}

\noindent
\textbf{Conjecture:} \cite{S} $\chi^A(M,N) \geq 0$.

\vspace{0.3cm}

\noindent
When $A$ is of finite type over a field, the above conjecture was answered in \cite{S}. Later a geometric proof was given by O. Gaber \cite{GA} for the non-negativity. One can express the Euler characteristic in terms of projective resolutions. If $E_{\bullet}$ and $F_{\bullet}$ be free resolutions of the $A$-modules $M,N$ (which may be taken to be finite, by the regularity of $A$), then

\begin{equation}
\chi^A(M,N) =  \chi (E_{\bullet} \otimes F_{\bullet})
\end{equation}

\noindent
where the right hand side is the usual Euler characteristic of the complex $E_{\bullet} \otimes F_{\bullet}$. The latter makes sense for the complex is supported on the maximal ideal of $A$. P. Roberts \cite{R} in his proof of the vanishing part of the conjecture uses the relation 

\begin{equation}
\chi(E_{\bullet})=\text{ch}(E_{\bullet}). \text{td}(A)
\end{equation}

\noindent
between the Euler characteristic and the local chern character. Then by the Riemann-Roch theorem

\begin{equation}
\chi( E_{\bullet} \otimes F_{\bullet})=ch(E_{\bullet} \otimes F_{\bullet})[A]
\end{equation}

\noindent
We try to explain this identity by a modification of the chern character 

\begin{equation}
ch:K_0(X) \to CH^*(X) \to H^*(X)
\end{equation}

\noindent
making the integrand a simple Poincare-Hodge dual pairing. The strategy goes through the definition of Mukai vector which we proceed to define.

\vspace{0.3cm}

\section{Mukai vector}

\vspace{0.3cm}

The references for this section are \cite{C}, \cite{HJLM} and \cite{F}. Let $X,Y$ be complex manifolds, and let $\Gamma \in D^b(X \times Y)$. Let $pr_1,pr_2$ be the projections. Define the Fourier transform with kernel $\Gamma$ by;

\begin{equation}
\Gamma_*: D^b(X) \to D^b(Y), \qquad \Gamma_*(.)=pr_{1,*}(pr_2^*(.) \otimes \Gamma)
\end{equation}

\vspace{0.2cm}

\noindent
Similarly one can define the Fourier transform on the cohomologies with kernels in $H^*(X \times Y,\mathbb{C})$. In order to relate them we use chern character and Riemann-Roch theorem. The Riemann-Roch theorem states that, associated to $\pi:X \to Y$ a smooth morphism;

\begin{equation}
\pi_*(\text{ch}(\bullet)\text{td}(X))=\text{ch}(\pi_*(\bullet)).\text{td}(Y)
\end{equation}

\noindent
where $ch:D^b(X) \to H^*(X)$ is the local chern character. We define the Mukai vector as follows,

\begin{equation}
\tilde{ch}:D^b(X) \to H^*(X,\mathbb{Q}), \qquad \tilde{ch}(.)=\text{ch}(.).\sqrt{\text{td}(X)}
\end{equation}

\noindent
where $\text{td}_X^{1/2} \in \oplus_{i} H^i(X,\Omega_X^i)$. The map in (8) is compatible with the Fourier transform defined in (6),

\begin{equation}
\begin{CD}
D^b(X) @>{\Gamma_*}>> D^b(Y) \\
@V{\tilde{ch}}VV                   @VV\tilde{ch}V\\
H^*(X,\mathbb{Q}) @>{\tilde{ch}(\Gamma)_*}>> H^*(Y,\mathbb{C})
\end{CD}
\end{equation}

\noindent
The map $\tilde{ch}(\Gamma)_*$ does respects the columns of Hodge diamond;

\begin{equation}
\tilde{ch}(\Gamma)_*: \displaystyle{\bigoplus_{p-q=i}H^{p,q}(X) \to \bigoplus_{p-q=i}H^{p,q}(X)}
\end{equation}

\noindent
That is because the class $\tilde{ch}(.)$ is a Hodge class. We define an involution on $H^*(X,\mathbb{C})$ as,

\begin{align}
\theta:\bigoplus_i H^i(X, \mathbb{C}) \longrightarrow \bigoplus_iH^i(X, \mathbb{C}) \qquad \\
\theta(v_0,v_1,...,v_{2n})=(v_0,iv_1,-v_2,...,i^{2n}v_{2n})
\end{align}

\noindent
It induces an operator

\begin{equation}
.^{\vee}:H^*(X,\mathbb{C}) \to H^*(X,\mathbb{C}), \qquad v^{\vee}=\theta(v).\displaystyle{\frac{1}{\sqrt{ch(\omega_X)}}}
\end{equation}

\noindent
One has $\text{td}(T_X^{\vee})=\text{td}(T_X).\exp(-c_1(T_X))=\text{td}(T_X).\text{ch}(\omega_X)$, where $\omega_X$ is the canonical sheaf of $X$. We have 

\begin{proposition} \cite{C} If ${E}$ and ${F}$ are coherent sheaves on the smooth projective variety $X$, 

\begin{equation}  
\chi({E},{F})=\langle \tilde{ch}(E), \tilde{ch}^{\vee}(F) \rangle 
\end{equation}
\end{proposition}

\noindent
We make a modification on the Mukai vector as follows. We replace the Mukai vector by the vector

\begin{equation} 
E \mapsto (2 \pi i)^{\deg(.)/2}\displaystyle{\frac{1}{(2 \pi)^{d/2}}}\Gamma(X) \wedge ch(E) 
\end{equation}

\noindent
The cohomology class ${\Gamma}_X$ is defined via the identity $\frac{z}{1-e^{-z}}=e^{i \pi z}\Gamma(1-x)\Gamma(1+x)$ used to share the two factors of $\sqrt{td_X}$ with the other chern classes in the Mukai pairing. It explicitly is given by the formula,

\begin{equation}
{\Gamma}_X=\exp(C.ch_1(T_X)+\displaystyle{\sum_{n \geq 2}\dfrac{\zeta(n)}{n}ch_n(T_X)}) 
\end{equation}

\noindent
where $C=\lim_{n \to \infty} (1+\frac{1}{2}+...+\frac{1}{n}-ln(n))$ is the Euler constant, $\zeta$ is the Riemann zeta. Let's write,

\begin{equation}
(2 \pi i)^{\deg(.)/2}{\hat{\Gamma}_X \wedge (\bullet)}: H^*(X,\mathbb{C}) 
{ \ {\rightarrow} H^*(X,\mathbb{C})}
\end{equation}

\noindent
We modify the vector $\nu({E})$ more and define 

\begin{equation}
\mu_{\Lambda}({E}):=ch({E}) \sqrt{td_X} . \exp(i\Lambda)
\end{equation}

\noindent
where $\Lambda$ is chosen so that $\theta(\Lambda)=-\Lambda$. The former Mukai vector is the special case $\Lambda=0$. We have 

\begin{equation} 
\chi({E},{F})=\displaystyle{\int ch({E}) \wedge ch({F})^{\vee} . td_X=\langle \mu_{\Lambda}({E}), \mu_{\Lambda}({F}) \rangle =\int \mu_{\Lambda}({E})\wedge \mu_{\Lambda}({F})^{\vee}} 
\end{equation}

\noindent
When $Y$ and $Z$ are projective sub-varieties of $X$, then the sheaves $\mathcal{E}, \mathcal{F}$ will be replaced by $\mathcal{O}_Y,\mathcal{O}_Z$, respectively. We denote the Euler characteristic by $\chi(Y,Z)$ in this case, for short.

\vspace{0.3cm}

\section{Multiplicity question over l-adic fields}

\vspace{0.3cm}

The references for this section are \cite{G} and \cite{SA}. Let $X$ be a smooth projective variety $/\mathbb{Q}_l$ of dimension $n$, and $L$ an ample divisor class. $L$ acts on etale cohomology of $X$ and by hard Lefschetz theorem, 

\begin{equation}
L^j:H^{n-j}(X(\bar{\mathbb{Q}_l}), \mathbb{Q}_l) \cong H^{n+j}(X(\bar{\mathbb{Q}_l}), \mathbb{Q}_l)
\end{equation}

\noindent
It follows that

\begin{equation}
H^{n-j}(X(\bar{\mathbb{Q}_l}), \mathbb{Q}_l) = \oplus_k L^k H^{j-2k}(X(\bar{\mathbb{Q}_l}), \mathbb{Q}_l)_{prim}
\end{equation}

\noindent
namely Lefschetz decomposition. The Grothendieck standard Conjecture $B$, says the Lefschetz decomposition (20) is in fact algebraic and is given over the Chow groupsby an algebraic cycle which we also denote by $L$, i.e.

\begin{equation}
A^j(X)=\oplus_k L^k.A^{j-k}(X)^{prim}
\end{equation}

\noindent
Furthermore, the pairing

\begin{equation}
(-1)^j\langle L^{n-2j}a,b \rangle , \qquad a,b \in A^j(X)^{prim}, \ \ j \leq n/2
\end{equation}

\noindent
There exists similar pairings on the cohomologies in (20) which are known to be positive definite by Hodge theory. They can also be defined by the Hodge Star operator 

\begin{equation}
*(L^km)=(-1)^{i(i+1)/2}L^{n-i-k}m 
\end{equation}

\noindent
on $H^*(X,\mathbb{Q}_l)$ by the pairing $(m,*n)$ and are compatible with the above identities. On the other hand, for a non-zero correspondence $\lambda \in A^n(X \times_{\mathbb{Q}_l} X) \subset End(H^*(X,\mathbb{Q}_l))$ we consider the transpose $\lambda^{\dagger}$ w.r.t this pairing. Then assuming the Grothendieck standard conjectures are satisfied, we will always have 

\begin{equation}
Tr(\lambda^{\dagger} \circ \lambda) \geq 0  
\end{equation} 

\noindent
This is analogue of positivity for Rosati involution, cf. \cite{SA} page 15. The action of the correspondences always determine the pairing on the homologies, by composing the action of diagonal $\Delta$ with the product structure. That is if $Y,Z \subset X$ be closed subvarieties 

\begin{equation}
\chi(Y,Z)=\chi(\Delta, Y \times Z)=\langle \Delta^{\vee}, Y \times Z \rangle =Tr (\Delta^{\dagger} \circ (Y \times Z))
\end{equation}

\noindent
Then one may be able to discuss the vanishing and positivity of $\chi(Y,Z)$ using the pairings (23) and (24).

\vspace{0.3cm}

\section{Appendix: Grothendieck Standard Conjectures}

\vspace{0.3cm}

We list the Grothendieck Standard conjectures, \cite{G}:

\begin{itemize}

\item A : Hard Lefschetz on cycles

\begin{equation}
A(X): .L^{n-2k}:CH^r(X) \cong CH^{n-r}(X)
\end{equation}

\item B : Lefschetz type Standard Conjecture

\begin{equation}
B(X): *L:\oplus_{i,r} H^i(X)(r) \to \oplus_{i,r}H^i(X)(r) \qquad is \ algebraic.
\end{equation}

\item C : Kunneth type Standard Conjecture

\begin{equation}
C(X): \pi_X^i:H^{\bullet}(X) \twoheadrightarrow H^{i}(X) \hookrightarrow H^{\bullet}(X) \qquad is \ algebraic
\end{equation}

\item D : Homological and numerical equivalence coincide

\begin{equation}
D(X): \qquad \sim_{hom,\mathbb{Q}} = \sim_{num,\mathbb{Q}}
\end{equation}

\item I : Hodge type Standard conjecture; the $\mathbb{Q}$-valued quadratic form $\alpha \mapsto \langle \alpha, *L(\alpha) \rangle $ on $Z_{hom}^{\bullet}(X)_{\mathbb{Q}}$ is positive definite.
\end{itemize}


\begin{thebibliography}{99}

\bibitem[C]{C}  A. Caldararu, The Mukai pairing and integral transforms in Hochschild homology, MOSCOW MATHEMATICAL JOURNAL
Volume 10, Number 3, July–September 2010, Pages 629–645

\bibitem[F]{F} W. Fulton, Intersection theory, Ergebnisse der Mathematik und ihrer Grenzgebiete. 3. Folge. A, Springer-Verlag, Berlin, 1998

\bibitem[GA]{GA} Gabber, O., Non-negativity of serre's intersection multiplicities, Expose a L'IHES (1995)

\bibitem[G]{G}  A. Grothendieck,  Standard conjectures on algebraic cycles, IHES, France, 1968

\bibitem[HJLM]{HJLM} J. Halverson, H. Jockers, J. Lapan, D. Morrison, Perturbative corrections to Kahler moduli spaces, UCSB Math 2013

\bibitem[R]{R}  P. Roberts, The vanishing of intersection multiplicities of perfect complexes, Bull. Amer. Math. Soc. (N.S.)
Volume 13, Number 2 (1985), 127-130

\bibitem[S]{S} J. P. Serre, Local Algebra, Multiplicities, Lecture notes in mathematics 11, Springer Verlag, Newyork, 1961

\bibitem[SA]{SA}  M. Saito, Monodromy filtration and positivity, RIMS, Kyoto University, arxiv:math/000162v6, 2000


\end{thebibliography}
\end{document}